\documentclass{article}
\usepackage{amssymb}
\usepackage{amscd}
\newtheorem{theorem}{Theorem}[section]
\newtheorem{remark}{Remark}[section]
\newtheorem{corollary}{Corollary}[section]
\newtheorem{definition}{Definition}[section]

\begin{document}

\begin{center}
\vskip 1cm{\LARGE\bf On a Class of Polynomials with Integer
Coefficients } \vskip 1cm \large {
Milan Janji\'c \\
Department of Mathematics and Informatics\\
 University of Banja Luka \\
Republika Srpska, Bosnia and Herzegovina} \\
email:agnus@blic.net
\end{center}

\vskip .2 in
\begin{abstract}
A class $P_{n,m,p}(x)$ of polynomials is defined. The
combinatorial meaning of its coefficients is given. Chebyshev
polynomials are the special cases of $P_{n,m,p}(x).$ It is first
shown that $P_{n,m,p}(x)$ may be expressed in terms of
$P_{n,0,p}(x).$ From this we derive that $P_{n,2,2}(x)$ may be
obtain in terms of trigonometric functions, from which we obtain
some of its important properties.

Some questions about orthogonality are also  concerned.

  Furthermore, it is shown that $P_{n,2,2}(x)$ fulfills the same three terms recurrence
as Chebyshev polynomials. Some others recurrences for
$P_{n,m,p}(x)$  and its coefficients are also obtained.

 At the end a formula for
coefficients of Chebyshev polynomials of the second kind is
derived.
\end{abstract}

\section{Introduction}
In the paper [1] the following result is proved.

\noindent\textbf{Theorem A.}\textit{ If a finite set $X$ consists
of $n$ blocks of the size $p$ and an additional block of the size
$m$ then, for $n\geq 0,\;k\geq 0,$ the number $f(n,k,m,p)$ of
$n+k-$ subsets of $X$ intersecting each block of the size $p$ is}
$$f(n,k,m,p)=\sum_{i=0}^n(-1)^i{n\choose i}{np+m-ip\choose n+k}.$$

The following  relations for the function $f$ are also proved in
[1]:
\begin{equation}\label{e1}f(n,k,m,p)=\sum_{i=0}^{m}{m\choose i}f(n,k-i,0,p),\end{equation}
\begin{equation}\label{e2}f(n,k,m,p)=\sum_{i=0}^t(-1)^i{t\choose i}f(n,k+t,m+t-i,p),\end{equation}
\begin{equation}\label{e3}f(n,k,m,p)=\sum_{i=1}^{p}{p\choose i}f(n-1,k-i+1,m,p-1),\end{equation}
\begin{equation}\label{e4}f(n,k,m,p)=\sum_{i=0}^{n}\sum_{j=0}^i{n\choose i}{i\choose j}f(n-j,k-i+j,m,p-1).\end{equation}

Furthermore, it is shown that  $(-1)^kf(n,k,0,2)$ is the
coefficient of Chebyshev polynomial $U_{n+k}(x)$ by $x^{n-k},$ and
that $(-1)^kf(n,k,1,2)$ is the coefficient of Chebyshev polynomial
$T_{n+k-1}(x)$ by $x^{n-k+1}.$

\begin{definition}
We define the set of coefficients
  $$\{c(n,k,m,p):\;n=m,m+1,\ldots;k=0,1,\ldots,n\}$$
 such that
$$c(n,k,m,p)=(-1)^{\frac{n-k}{2}}f\left(\frac{n+k-2m}{2},\frac{n-k}{2},m,p\right),$$
if $n$ and $k$ are of the same parity, and $c(n,k,m,p)=0$
otherwise. Polynomials $P_{n,m,p}(x)$ are defined to be
$$P_{n,m,p}(x)=\sum_{k=0}^nc(n,k,m,p)x^k.$$
\end{definition}

\begin{remark}
 Chebyshev polynomials are particular cases of $P_{n,m,p}(x),$ obtained for $m=1,\;p=2$ and $m=0,\;p=2,$
 that is,
 $$U_n(x)=P_{n,0,2},\;T_n(x)=P_{n,1,2}.$$
 \end{remark}

The polynomial $P_{n,2,2}(x)$ is the closest to Chebyshev
polynomials, and will be denoted simply by $P_n(x).$

  In the next table we state the first few of $P_{n}(x).$
$$\begin{array}{c}x^2\\2x^3-2x\\4x^4-5x^2+1\\8x^5-12x^3+4x\\16x^6-28x^4+13x-1\\
32x^7-64x^5+38x^3-6x.\end{array}$$

Among coefficients of the above polynomials the following
sequences from [2] appear: A024623, A049611, A055585, A001844,
A035597.

Triangles of coefficients for $P_{n,m,2}(x),(\;m=2,3,4,5,6)$ are
given in  A136388, A136389, A136390, A136397, A136398
respectively.

\section{Reduction to the case $m=0.$}

We shall first prove an  analog of the formula (\ref{e1})  for
polynomials.

\begin{theorem} The following equation is fulfilled:
$$P_{n,m,p}(x)=\sum_{i=0}^m(-1)^{i}{m\choose i}x^{m-i}P_{n-m-i,0,p}(x).$$
\end{theorem}

\noindent\textbf{Proof.} It holds
$$P_{n,m,p}(x)=\sum_{k=0}^n(-1)^{\frac{n-k}{2}}
f\left(\frac{n+k-2m}{2},\frac{n-k}{2},m,p\right)x^k.$$ Using
(\ref{e1}) one obtains
$$P_{n,m,p}(x)=\sum_{k=0}^n\sum_{i=0}^{m}(-1)^{\frac{n-k}{2}}{m\choose i}f\left(r,s,0,p\right)x^k.$$
where $$r=\frac{n+k-2m}{2},\;s=\frac{n-k}{2}-i.$$ Changing the
order of summation yields
$$ P_{n,m,p}(x)=\sum_{i=0}^{m}{m\choose
i}x^{m-i}\sum_{k=0}^n(-1)^{\frac{n-k}{2}}f\left(r,s,0,p\right)x^{k-m+i}.$$

Terms in the sum on the right side of the preceding equation
produce nonzero coefficients only in the case $0\leq s\leq r,$
that is,
$$m-i\leq k\leq n-2i.$$ It follows that

$$P_{n,m,p}(x)=\sum_{i=0}^{m}(-1)^i{m\choose
i}x^{m-i}\sum_{k=m-i}^{n-2i}(-1)^{\frac{n-k}{2}-i}f\left(r,s,0,p\right)x^{k-m+i}.$$

Denoting $k-m+i=j$ we obtain
$$P_{n,m,p}(x)=\sum_{i=0}^{m}(-1)^{i}{m\choose
i}x^{m-i}\sum_{j=0}^{n-i}c(n-i,j,0,p)x^{j},$$ which means that

$$P_{n,m,p}(x)=
\sum_{i=0}^{m}(-1)^{i}{m\choose i}x^{m-i}P(n-m-i,0,p)(x),$$ and
the theorem is proved.

According the the preceding theorem we may express $P_n(x)$ in
terms of Chebyshev polynomials of the second kind.  Namely, for
$m=2,\;n\geq 4$ holds
\begin{equation}\label{e5}P_{n}(x)=x^2U_{n-2}(x)-2xU_{n-3}(x)+U_{n-4}(x).\end{equation}
This allow us to express $P_{n}(x)$ in terms of trigonometric
functions.

\begin{theorem}
For each  $n\geq 3$ holds
\begin{equation}\label{ee}P_{n}(\cos\theta)=-\sin\theta\sin(n-1)\theta.\end{equation}
\end{theorem}

\noindent\textbf{Proof.} According to (\ref{e5}) and well-known
property of  Chebyshev polynomials we obtain
 $$\sin\theta P_{n}(\cos\theta)=\cos^2\theta\sin(n-1)\theta-2\cos\theta\sin(n-2)\theta+\sin(n-3)\theta.$$
From the identity
$$2\cos\theta\sin(n-2)\theta=\sin(n-1)\theta+\sin(n-3)\theta$$
follows
$$\sin\theta P_{n}(\cos\theta)=\cos^2\theta\sin(n-1)\theta-\sin(n-1)\theta=-\sin^2\theta\sin(n-1)\theta.$$
Dividing by $\sin\theta\not=0$ we prove the theorem.

Note that this proof is valid for $n\geq 4.$ The case $n=3$ may be
checked directly.

In the following theorem we prove that $P_{n}(x)$ have  the same
important property concerning zeroes as Chebyshev polynomials do.

\begin{theorem}
For $n\geq 3,$ the polynomial $P_{n}(x)$ has all simple zeroes
lying in the segment $[-1,1].$
\end{theorem}

\noindent\textbf{Proof.} Since
$$U_n(1)=n+1,\;U_n(-1)=(-1)^n(n+1)$$ the equation (\ref{e5})
implies

$$P_{n}(1)=U_{n-2}(1)-2U_{n-3}(1)+U_{n-4}(1)=n-1-2(n-2)+n-3=0,$$
and
$$P_{n}(-1)=U_{n-2}(-1)+2U_{n-3}(-1)+U_{n-4}(-1)=$$$$=(-1)^{n-2}(n-1)+2(-1)^{n-3}(n-2)+(-1)^{n-4}(n-3)=0.$$

Thus, $x=-1$ and $x=1$ are zeroes of $P_{n}(x).$ The remaining
$n-2$ zeroes are obtained from the equation
$$\sin(n-1)\theta=0,$$ and they are
$$x_k=\cos\frac{k\pi}{n-1},\;(k=1,2,\ldots,n-2).$$

We shall now state an immediate consequence of (\ref{ee}) which
shows that  values of $P_{n}(x),\;(x\in[-1,1])$ lie inside the
unit circle.

\begin{corollary} For $n\geq 3$ and $x\in[-1,1]$ we have
$$P_{n}(x)^2+x^2\leq 1.$$
\end{corollary}

\begin{remark} Dividing $P_n(x)$ by $2^{n-2}$ we obtain a polynomial with the leading coefficient $1.$
  Thus, its  supremum  norm on $[-1,1]$ is
$\leq\frac{1}{2^{n-2}},$ which  means that
$\frac{1}{2^{n-2}}P_{n}(x)$ has at most $2$ times greater supremum
norm, comparing with the supremum norm of $T_{n}(x),$ that is
minimal.
\end{remark}

Taking derivative in the equation (\ref{ee}) we obtain the
following  equation for extreme points of $P_{n}(x):$

 $$(n-1)\tan\theta+\tan(n-1)\theta=0.$$

The values  $\theta=0,$ and $\theta=\pi$ obviously satisfied this
equation,
 which implies that endpoints $x=-1$ and $x=1$ are extreme points.
 The remaining extreme points of $P_3(x)$ are $x=\arctan\sqrt 2$ and
 $x=-\arctan\sqrt 2.$

\section{Orthogonality}

In this section we investigate the set
$\{P_{n}(x):n=2,3,4,\ldots\}$ concerning to the problem of
orthogonality, with respect to some standard Jacobi's weights.

The first result is for the weight $\frac{1}{\sqrt{1-x^2}}$ of
Chebyshev polynomials of the first kind.

\begin{theorem} It holds
$$\int_{-1}^1\frac{P_n(x)P_m(x)}{\sqrt{1-x^2}}dx=\left\{\begin{array}{rl}\frac{\pi}{4}&m=n\\
-\frac{\pi}{8}&|n-m|=2\\0&\mbox{ otherwise}.\end{array}\right.$$
\end{theorem}

\noindent\textbf{Proof.} Puting $x=\cos\theta$ implies
$$I=\int_{-1}^1\frac{P_n(x)P_m(x)}{\sqrt{1-x^2}}dx=\int_{0}^{\pi} P_n(\cos\theta)P_m(\cos\theta)d\theta.$$
Using (\ref{ee}) we
obtain$$I=\int_{0}^{\pi}\sin^2\theta\sin(n-1)\theta\sin(m-1)\theta
d\theta.$$ Transforming the integrating function we obtain
$$\sin^2\theta\sin(n-1)\theta\sin(m-1)\theta=
\frac14\cos(n-m)\theta-\frac14 \cos(n+m-2)\theta-$$$$-\frac18
\cos(n-m-2)\theta-\frac 18 \cos(n-m+2)\theta+
\frac18\cos(n+m-4)\theta+\frac18\cos(n+m)\theta.$$
 Taking into account that $m,n\geq 3$ we conclude that integrals of the terms on the right side of the preceding equation
 are zero if $n\not=m$ and
$|n-m|\not=2.$ If $n=m$ we obtain $I=\frac{\pi}{4},$ and
$I=-\frac{\pi}{8}$ if $|n-m|=2,$ and the theorem is proved.

\begin{corollary}Each subset of the set $\{P_n(x):n\geq 3\},$ not containing
polynomials $P_k(x)$ and $P_m(x)$ such that $|k-m|=2,$ is
orthogonal.
\end{corollary}

The next result concerns the weight $\sqrt{1-x^2}$ of Chebyshev
polynomials of the second kind. The result is similar to the
result of  the preceding theorem.

\begin{theorem}It holds
$$\int_{-1}^1\sqrt{1-x^2}P_n(x)P_m(x)dx=\left\{\begin{array}{rl}\frac{3\pi}{16}&m=n\\
-\frac{\pi}{8}&|n-m|=2\\\frac{\pi}{32}&|n-m|=4\\0&\mbox{
otherwise}.\end{array}\right.$$
\end{theorem}
\noindent\textbf{Proof.} In this case we have
$$\int_{-1}^1\sqrt{1-x^2}P_n(x)P_m(x)dx=\int_{0}^{\pi} \sin^2\theta P_n(\cos\theta)P_m(\cos\theta)d\theta.$$
We therefore  need to calculate the integral
$$\int_{0}^{\pi} \sin^4\theta\sin(n-1)\theta\sin(m-1)\theta d\theta.$$

In this case we have
$$\sin^4\theta\sin(n-1)\theta\sin(m-1)\theta=\frac{3}{16}\cos(n-m)\theta-\frac{3}{16}\cos(n+m-2)\theta+$$$$+
\frac{1}{32}\cos(n-m-4)\theta+\frac{1}{32}\cos(n-m+4)\theta-\frac{1}{32}\cos(n+m-6)\theta-\frac{1}{32}\cos(n+m+2)\theta-
$$$$-\frac{1}{8}\cos(n-m-2)\theta-\frac{1}{8}\cos(n-m+2)\theta+\frac{1}{8}\cos(n+m-4)\theta+\frac{1}{8}\cos(n+m)\theta.$$

The integral of each term on the right side with
$m\not=n,\;|m-n|\not=2,\;|n-m|\not=4$ is zero.

For these particular values we easily obtain the desired result,
and the theorem is proved.

Taking, for instance, the weight $(1-x^2)^{\frac 32}$ in the
similar way one obtains

$$\int_{-1}^1(1-x^2)^{\frac32}P_n(x)P_m(x)dx=\left\{\begin{array}{rl}\frac{5\pi}{32}&m=n\\
-\frac{15\pi}{128}&|n-m|=2\\\frac{3\pi}{64}&|n-m|=4\\-\frac{\pi}{128}&|n-m|=6\\0&\mbox{
otherwise}.\end{array}\right.$$

Considering the weight $1$ leads to the following result:
\begin{theorem} If $m$ and $n$ are of different parity then
$$\int_{-1}^1P_n(x)P_m(x)dx=0.$$
\end{theorem}

\noindent\textbf{Proof.} In this case we need to calculate the
integral
$$\int_{0}^{\pi} \sin^3\theta\sin(n-1)\theta\sin(m-1)\theta d\theta.$$
 We have
$$\sin^3\theta\sin(n-1)\theta\sin(m-1)\theta=
-\frac{1}{16}\sin(n-m+3)\theta+\frac{1}{16}\sin(n-m-3)\theta+$$$$+\frac{1}{16}\sin(n+m+1)\theta-\frac{1}{16}
\sin(n+m-5)\theta+\frac{3}{16}\sin(n-m+1)\theta-$$$$-\frac{3}{16}\sin(n-m-1)\theta-\frac{3}{16}\sin(n+m-1)\theta
+\frac{1}{16}\sin(n+m-3)\theta.$$ Since $m$ and $n$ are of
different parity each function on the right is of the form
$\sin(2k+1)\theta,$ which implies that  its integral is zero, and
the theorem is proved.

\section{Some recurrence relations }

In this section we prove some recurrence relation for
$P_{n,m,p}(x)$ as well as some recurrence relations for their
coefficients.
\begin{theorem} For each integer $t\geq 0$ holds
$$P_{n,m,p}(x)=
\sum_{i=0}^{t}(-1)^{t-i}{t\choose i}x^{i}P_{n+2t-i,m+t-i,p}(x).$$
\end{theorem}

\noindent\textbf{Proof.} Translating  (\ref{e2}) into the equation
for coefficients we obtain

$$c(n,k,m,p)=\sum_{i=0}^t(-1)^{i+t}{t\choose i}c(n+2t-i,k-i,m+t-i,p).$$

Multiplying by $x^k$ yields
$$c(n,k,m,p)x^{k}=\sum_{i=0}^t(-1)^{i+t}{t\choose i}x^ic(n+2t-i,k-i,m+t-i,p)x^{k-i},$$
which easily implies the claim of the theorem.

In the  case $t=1,\;m=1,\;p=2$ we obtain the following formula,
expressing $P_n(x)$ in terms of Chebyshev polynomials of the first
kind:
$$P_{n}(x)=xT_{n-1}(x)-T_{n-2}(x).$$

From this we easily conclude that $P_n(x)$ satisfies the same
three term recurrence as Chebyshev polynomials.

\begin{corollary} The polynomials $P_{n,m,2}(x)$ satisfy the following
equation:
$$P_{n,m,2}(x)=2xP_{n-1,m,2}(x)-P_{n-2,m,2}(x),$$
with initial conditions
 $$P_{0,m,2}(x)=x^m,\;P_{1,m,2}(x)=2x^{m+1}-mx^{m-1}.$$
\end{corollary}

Combining the equations (\ref{e1}) and (\ref{e4}) we obtain
$$f(n,k,m,p)=\sum_{i=0}^{n}\sum_{j=0}^i\sum_{t=0}^m{n\choose i}{m\choose t}{i\choose j}f(n-j,k-i+j-t,0,p-1).$$

Translating this equation into the equation for coefficient we
obtain

$$c(n,k,m,p)=\sum_{i=0}^{n}\sum_{j=0}^i\sum_{t=0}^m{n\choose i}{m\choose t}{i\choose j}(-1)^{i-j+t}c(n-i-t,k+i-2j+t,0,p-1).$$

Applying the preceding equation several times we obtain the
following:
\begin{corollary}
 Coefficients of $P_{n,m,p}(x)$
may be obtained as a functions of coefficients of Chebyshev
polynomials of the second kind.
\end{corollary}

Converting (\ref{e3}) into the equation for coefficients we obtain
$$c(n,k,m,p)=\sum_{i=1}^p(-1)^i{p\choose i}c_{n-i,k+i-2,m,p}.$$
This implies the following:
\begin{corollary}
Coefficients of $P_{n,m,p}(x)$ may be expressed in terms of
coefficients of polynomials $P_{n',m,p}(x),$ where $n'<n.$
\end{corollary}

We  shall finish the paper with a formula for coefficients of
Chebyshev polynomials of the second kind. Taking $p=2$ in
(\ref{e4}) we obtain
$$f(n,k,m,2)=\sum_{i=0}^{n}\sum_{j=0}^i{n\choose i}{i\choose j}f(n-j,k-i+j,m,1).$$

Since $f(r,s,m,1)={m\choose s}$ we have

$$f(n,k,m,2)=\sum_{i=0}^n\sum_{j=0}^i{n\choose i}{i\choose j}{m\choose k-i+j}.$$
For  $m=0,$ in the sum on the right side of this equation  only
terms with $k=i-j$ remains. We thus obtain
$$f(n,k,0,2)=\sum_{s=0}^{n-k}{n\choose s}{n-s\choose k}.$$
Accordingly, the following formula follows

\begin{corollary} For coefficients $c(n,k)$ of Chebyshev
polynomial $U_n(x)$ hold
$$c(n,k)=(-1)^{\frac{n-k}{2}}\sum_{i=0}^k{\frac{n+k}{2}\choose i}
{\frac{n+k}{2}-i\choose \frac{n-k}{2}},$$ if  $n$ i $k$ are of the
same parity and $c(n,k)=0$ otherwise.
\end{corollary}

\end{document}